\documentclass[a4paper,12pt]{article}
\usepackage[english]{babel}
\usepackage[T2A]{fontenc}
\usepackage[cp1251]{inputenc}
\usepackage[english]{babel}
\usepackage{amsthm}
\usepackage[tbtags]{amsmath}
\usepackage{amsfonts,amssymb}
\begin{document}

\newtheorem{theorem}{Theorem}
\newtheorem{lemma}{Lemma}
\newtheorem{proposition}{Proposition}
\newtheorem{Cor}{Corollary}

\begin{center}
{\bf On Realization and Isomorphism Problems\\ for Formal Matrix Rings}
\end{center}
\hfill \textbf{Piotr Krylov}, 

\hfill Tomsk National Recearch $\qquad\;\,$

\hfill State University, krylov@math.tsu.ru

\hfill \textbf{Askar Tuganbaev}

\hfill National Recearch University <<MPEI>>, 

\hfill Lomonosov Moscow State University 

\hfill tuganbaev@gmail.com

\textbf{Abstract.} We consider realization and isomorphism problems for formal matrix rings over a given ring. Principal multiplier matrices of such rings play an important role in this case.\\
The work of  A.A.Tuganbaev is supported by Russian Scientific Foundation, project 22-11-00052.

\textbf{Key words:} formal matrix ring, principal multiplier matrix

\textbf{MSC Classification. 16R99; 16D10}


\section{Introduction}\label{section1} 

Formal matrix rings  (or generalized matrix rings) over a given ring attract a lot of attention from specialists. This is natural, since such rings regularly appear in the theory of rings and modules. In particular, they play an important role in the study of a number of classes of Artinian rings and algebras (see \cite{AusRS95}, \cite{BabO09}). They also provide a variety of examples for the general theory of rings and modules. A number of aspects of the theory of formal matrix rings are presented in the book \cite{KryT17}.

There is one interesting form of formal matrix rings. In the case of $2\times 2$ matrices, they appeared in \cite{Kry08}; in the case of $n\times n$ matrices, they appeared in \cite{TanZ13}. We mean formal matrix rings over a given ring $R$ (or one says <<with values ??in $R$>>). This means that a particular formal matrix ring contains the same ring $R$ in all positions. The class of such rings is a direct extension of the usual ring of $n\times n$ matrices $M(n,R)$. However, the properties of formal matrix rings over the ring $R$ may be very different from the properties of the ring $M(n,R)$. Chapter 3 of the book \cite{KryT17} contains an exposition of some questions about formal matrix rings over $R$, and at the beginning of this chapter three problems about formal matrix rings are formulated. In chapters 18-20, these problems are solved for some types of formal matrix rings over $R$. This paper is devoted to two of the three indicated tasks. Namely, the realization problem I and the isomorphism problem III. These problems are considered in\cite{AbyT15}, \cite{AbyT15b}, \cite{CheDS20}, \cite{Tap15}, \cite{Tap17}. 

In this paper, we consider only associative rings with non-zero identity element. If $R$ is a ring, then $M(n,R)$ is the ordinary ring of all $n\times n$ matrices with values in the ring $R$. For an arbitrary ring $S$, the prime radical of $S$ is denoted by $P(S)$.

\section{Formal Matrix Rings over a Given Ring}\label{section2} 

Let us briefly recall the definition of a formal matrix ring. We fix a positive integer $n\ge 2$. Let $R_1,\ldots,R_n$ be $n$ rings and let $M_{ij}$ be $R_i$-$R_j$-bimodules such that $M_{ii}=R_i$, $i,j=1,\ldots,n$. We assume that for any subscripts $i,j,k=1,\ldots,n$, we have an $R_i$-$R_k$-bimodule homomorphism $M_{ij}\otimes_{R_j}M_{jk}\to M_{ik}$. We denote by $K$ the set of all $n\times n$ matrices with values in the bimodules $M_{ij}$. The set $K$ is a ring with respect to standard matrix operations of addition and multiplication. Matrices are multiplied by applying the bimodule homomorphisms mentioned above. The ring $K$  is called a \textsf{formal matrix ring} (or a \textsf{generalized matrix ring}) of order $n$.  It can be written in the following form:
$$
K=\begin{pmatrix}
R_1&M_{12}&\ldots&M_{1n}\\
M_{21}&R_{2}&\ldots&M_{2n}\\
\ldots&\ldots&\ldots&\ldots\\
M_{n1}&M_{n2}&\ldots&R_{n}
\end{pmatrix}.
$$
Let $R$ be a ring. If $K$ is a formal matrix ring such that $M_{ij}=R$ for all $i$ and $j$, then it is called a \textsf{formal matrix ring over the ring} $R$ or \textsf{formal matrix ring with values in the ring} $R$. Such rings can be defined directly. Namely, let
$\{s_{ijk}\,|\,i,j,k=1,\ldots,n\}$ be some set of central elements of the ring $R$ such that 
$$
s_{iik}=1=s_{ikk},\; s_{ijk}\cdot s_{ik\ell}=s_{ij\ell}\cdot s_{jk\ell} \eqno (1)
$$ 
for all subscripts $i,j,k,\ell=1,\ldots,n$. For arbitrary $n\times n$ matrices $A=(a_{ij})$ and $B=(b_{ij})$ with values in $R$, we define a new multiplication by setting
$$
AB=C=(c_{ij}),\, \text{ where } c_{ij}=\sum_{k=1}^ns_{ikj}a_{ik}b_{kj}.
$$
As a result, we obtain a ring which is denoted by $K$ or $M(n,R,\Sigma)$, where $\Sigma$ is the set of all $s_{ijk}$. The set $\Sigma$ is called the \textsf{multiplier system} and the elements of $\Sigma$ are called \textsf{multipliers} of the ring $K$. If all $s_{ijk}$ are equal to $1$, then we obtain the ordinary matrix ring $M(n,R)$.

Let $\tau$ be a permutation of degree $n$. If $\Sigma=\{s_{ijk}\}$ is some multiplier system, then we set $t_{ijk}=s_{\tau(i)\tau(j)\tau(k)}$.  Then $\{t_{ijk}\}$ also is a multiplier system, since it satisfies identities $(1)$. We denote it by $\tau\Sigma$. Consequently, there exists a formal matrix ring $M(n,R,\tau\Sigma)$. The rings $M(n,R,\Sigma)$ and $M(n,R,\tau\Sigma)$ are isomorphic to each other under the correspondence $A\to \tau A$, where $A=(a_{ij})$ and $\tau A=(a_{\tau(i)\tau(j)})$.

Several matrices can be associated with the ring $M(n,R,\Sigma)$. We set $S=(s_{iji})$ and $S_k=(s_{ikj})$, for any $k=1,\ldots,n$. These matrices are called \textsf{multiplier matrices} of the ring $M(n,R,\Sigma)$. The matrix $S$ is
symmetrical. Following \cite{CheDS20}, we call it a \textsf{principal multiplier matrix}. In \cite{CheDS20}, the matrices $(s_{ijk})$ and $(s_{kij})$ are also used, $k=1,\ldots,n$. It is clear that  the matrices $\tau S$ and $\tau S_k$ are the corresponding multiplier matrices of the ring $M(n,R,\tau\Sigma)$.

\textbf{Until the end of the article, we assume that $K$ is a formal matrix ring over a given ring $R$. In addition, with the exception of Theorem 4.3, every multiplier $s_{ijk}$ is equal to $0$ or $1$. Therefore, all multiplier matrices of such a ring are $(01)$-matrices. For brevity, we call the ring $K$ a \textsf{$(01)$}-formal matrix ring.}

The following result can be easily obtained from identities $(1)$.

\textbf{Lemma 2.1.} For multipliers $s_{iji}$, $s_{jkj}$, and  $s_{kik}$ of some $(01)$-formal matrix ring, one of the following possibilities takes place (where $i\ne j$, $j\ne k$, and $k\ne i$).

\textbf{1)} All three elements are equal to $1$.

\textbf{2)} Some two of these three elements are equal to $0$ and the third element is $1$.

\textbf{3)} All three elements are equal to $0$.

It is useful to introduce the notion of an abstract multiplier matrix. Let $T=(t_{ij})$ be a symmetrical $n\times n$~ $(01)$-matrix over the ring $R$ such that all elements of the main diagonal are equal to $1$ and one of assertions 1)--3) of Lemma 2.1 is true for any three elements $t_{ij}$, $t_{jk}$, $t_{ki}$ with pairwise distinct subscripts. We call such a matrix $T$ a \textsf{principal multiplier matrix}.

\textbf{Lemma 2.2 \cite{KryN18}, \cite[Chapter 12]{KryT21}.} For any principal $n\times n$ multiplier matrix $T$, there exists a permutation $\sigma$ of degree $n$ such that the matrix $\sigma T$ is of canonical form.

By \textsf{canonical form} it is meant that the matrix $\sigma T$ can be represented in such a block form that blocks of 1's are on the main diagonal, and zeros are in all other positions.

Under the conditions of Lemma 2.2, we say that the matrix $T$ \textsf{is reduced to the canonical form}. It follows from Lemma 2.1 that the principal multiplier matrix of any $(01)$-formal matrix ring can be reduced to the canonical form.

\textbf{Remark 2.3.} We clarify that the canonical form is determined up to permutation of blocks on the main diagonal and up to the corresponding permutation of remaining blocks.

\section{Realization Problem for $(01)$-Formal Matrix Rings}\label{section3}

Section 1 contains information on the realization problem. The problem is related to description of multiplier matrices as abstract matrices.

\textbf{Theorem 3.1.} Let $T$ be some principal multiplier matrix. There exists a $(01)$-formal matrix ring such that principal multiplier matrix coincides with $T$.

\textbf{Proof.} Based on Lemma 2.2, we assume that the matrix $T$ is of canonical form. For all subscripts $i,j,k\in\{1,\ldots,n\}$, we define an element $s_{ijk}$ of the ring $R$ such that $s_{ijk}=1$ if one of the pairs $(i,j)$ or $(j,k)$ takes a position in some block on the main diagonal of the matrix $T$. Otherwise, we assume that $s_{ijk}=0$. The set $\Sigma$ of all such elements $s_{ijk}$ satisfies identities $(1)$. Therefore, the ring $M(n,R,\Sigma)$ exists. The matrix $T$ is a principal multiplier matrix of the ring $M(n,R,\Sigma)$.~$\square$ 

Multipliers of the ring $M(n,R,\Sigma)$, constructed in the proof of Theorem 3.1, satisfy the following condition:
\begin{itemize}
\item $s_{ijk}=0$ for any pairwise distinct subscripts $i,j,k$ such that $s_{iji}=0=s_{jkj}$.
\end{itemize}
We say that some $(01)$-formal matrix ring is a ring \textsf{of the form $K_0$} if its multipliers satisfy the above condition $\bullet$. By the use of a principal multiplier matrix of such a ring,  all remaining multiplier matrices can be  restored.

We note that papers \cite{KryN18}, \cite{KryT21} and \cite{KryT22} contain various material on automorphism groups of $(01)$-formal matrix rings.

\section{Isomorphism Problem for Formal Matrix Rings}\label{section4}

In \cite[Chapter 16]{KryT17}, the following isomorphism problem III for formal matrix rings is formulated:
\begin{itemize}
\item When two multiplier systems define isomorphic formal matrix rings?
\end{itemize}
We consider this problem for $(01)$-formal matrix rings.

We say that a ring $R$ \textsf{satisfies $(n,m)$-condition} if for any positive integers $n$ and $m$, we have $m=n$ provided  the rings $M(n,R)$ and $M(m,R)$ are isomorphic. For example, the $(n,m)$-condition holds if the ring $R$ is commutative, or local, or is a principal left (right) ideal domain.

A ring $S$ is said to be \textsf{indecomposable} if $0$ and $1$ are only central idempotents of $S$.

\textbf{Theorem 4.1.} We assume that the factor ring $R/P(R)$ is indecomposable and satisfies the $(n,m)$-condition and $K_1$, $K_2$ are $(01)$-formal matrix rings with principal multiplier matrices $S$ and $T$, respectively. The following assertions are true.

\textbf{1.} If the rings $K_1$ and $K_2$ are isomorphic to each other, then the matrices $S$ and $T$ have the same canonical forms.

\textbf{2.} If $K_1$ and $K_2$ are $(01)$-formal matrix rings of the form $K_0$ and the canonical forms of the matrices $S$ and $T$ coincide, then the rings $K_1$ and $K_2$ are isomorphic to each other.

\textbf{Proof.} \textbf{1.} By Lemma 2.2, we can assume that the matrices $S$ and $T$ are given in the canonical form. Let us assume that $K_1\cong K_2$. Then there exists a ring isomorphism
$$
\gamma\colon K_1/P(K_1)\to K_2/P(K_2).
$$
The structure of the prime radicals $P(K_1)$ and $P(K_2)$ is known (see \cite[Corollary 17.1]{KryT17} and the paragraph after the corollary). We also know the block structure of matrices $S$ and $T$. With the use of this information, we obtain the relations
$$
K_1/P(K_1)=P_1\times\ldots\times P_k\, \text{ and }\, K_2/P(K_2)=Q_1\times\ldots\times Q_{\ell},
$$
where $k$ (resp., $\ell$) the number of blocks on the main diagonal of the matrix $S$ (resp., $T$). In addition, all $P_i$ and $Q_j$ are full matrix rings of some orders. Since the ring $R/P(R)$ is indecomposable, all the rings $P_i/P(P_i)$ and $Q_j/P(Q_j)$ are indecomposable.

At this point, we note that there is an analogue of \cite[Lemma 9.6]{KryT21} (or \cite[Lemma 6.1]{KryT22}) on automorphisms of direct products of indecomposable rings for isomorphisms between direct products of indecomposable rings. Therefore, $k=\ell$ and there exists a permutation $\tau$ of degree $k$ such that the restriction $\gamma$ to $P_i$ is an isomorphism $P_i\to Q_{\tau(i)}$, $i=1,\ldots,k$. Consequently, canonical forms matrices $S$ and $T$ coincide.

\textbf{2.} Let $\{s_{ijk}\}$ (resp., $\{t_{ijk}\}$) be the set of all multipliers of the ring $K_1$ (resp., $K_2$). As noted, multipliers  of the form $s_{iji}$, i.e., elements of the principal multiplier matrix $S$, define all remaining multipliers $s_{ijk}$; the same is true for multipliers  of the ring $K_2$. Therefore, $s_{ijk}=t_{ijk}$, for all $i,j,k$. Therefore, we have the relation $K_1=K_2$.~$\square$

\textbf{Corollary 4.2.} The factor rings $K_1/P(K_1)$ and $K_2/P(K_2)$ are isomorphic to each other if and only if the matrices $S$ and $T$ have the same canonical form.

We introduce one more class of formal matrix rings over a given ring, we also prove some theorem related to the  isomorphism problem, for such rings.

We fix some central element $s$ of the ring $R$ such that $s^2\ne 1$ and  $s^2\ne s$. In addition, let $K$ be a formal matrix ring over $R$ such that every multiplier of $K$ is equal to $1$ or $s$. We denote the defined ring by $M(n,R,s)$ and call it a \textsf{$(s1)$-formal matrix ring}. The rings $M(n,R,s)$ are studied in Chapters 18--19 of the book \cite{KryT17}. Relationships between the rings $M(n,R,s)$ and crossed matrix ring are known (see \cite{BabO09}). In \cite[Lemma 4.7]{CheDS20}, it is proved that the principal multiplier of the matrix ring $M(n,R,s)$ determines all other multiplier matrices. For $(01)$-formal matrix rings, a similar assertion does not appear to be true.

Let $S$ be the principal multiplier matrix of some ring $M(n,R,s)$. According to \cite[Lemma 18.2]{KryT17}, there exists a permutation $\tau$ of degree $n$ such that the matrix $\tau S$ is of canonical form. This means that the matrix $\tau S$ can be presented in a block form such that the main diagonal contains blocks consisting of $1$s and the element $s$ is in all remaining positions.

\textbf{Theorem 4.3.} Let the factor ring $R/P(R)$ be indecomposable and satisfy the $(n,m)$-condition. For $(s1)$-formal matrix rings $K_1$ and $K_2$ over $R$, where $s\in P(R)$, the following assertion is true.\\
The rings $K_1$ and $K_2$ are isomorphic to each other if and only if for the rings $K_1$ and $K_2$, the canonical forms of principal multiplier matrices coincide.

\textbf{Proof.} We assume that principal multiplier matrices for the rings $K_1$ and $K_2$ already are of canonical form.

\textsf{Necessity.} If consider that $s\in P(R)$, then, in fact, we can repeat the proof of Theorem 4.1(1).

\textsf{Sufficiency.} Principal multiplier matrices for the rings $K_1$ and $K_2$ uniquely define all remaining multiplier matrices, see \cite[Lemma 4.7]{CheDS20}. Therefore, we can repeat the argument from the proof of Theorem 4.1(2).~$\square$

\textbf{Remark 4.4.} In \cite[Theorem 4.12]{CheDS20}, a result, which is similar to Theorem 4.3, is proved for a left Artinian ring $R$. 

\label{biblio}

\end{document}